\documentclass[11pt,a4paper]{article}
\usepackage[a4paper,
            bindingoffset=0.2in,
            left=1in,
            right=1in,
            top=1in,
            bottom=1in,
            footskip=.25in]{geometry}
\usepackage{amssymb}
\usepackage{cite}
\usepackage[hyphens]{url} 
\usepackage{amsthm}
\usepackage{epstopdf}
\usepackage{float}
\usepackage{graphicx}
\usepackage{amsmath}
\usepackage{authblk}

\theoremstyle{plain}
 \newtheorem{thm}{Theorem}[section]

\theoremstyle{definition}
 \newtheorem{exm}{Example}[section]
 
\theoremstyle{remark}
\newtheorem{rem}{Remark}[section]
\numberwithin{equation}{section}

\title{Converting Exhausters and Coexhausters}

\author[1]{M.~E. Abbasov\thanks{m.abbasov@spbu.ru}}
\affil[1]{St. Petersburg State University, SPbSU, 7/9 Universitetskaya nab., St. Petersburg, 199034 Russia}

\begin{document}

\maketitle

\begin{abstract}
We develop sufficient conditions for the equality of $\min\max$
and $\max\min$ in case when maximums and minimums are taken over
finite sets. We formulate these conditions in a form employable for the
construction of a new converting procedure for exhausters and
coexhausters.

Exhausters and coexhausters are notions of constructive nonsmooth
analysis which are used to study extremal properties of functions.
An upper exhauster (coexhauster) is used to get an approximation
of a considered function in the neighborhood of the point in the
form of $\min\max$ of linear (affine) functions. A lower exhauster
(coexhauster) is used to represent the approximation in the form
of $\max\min$ of linear (affine) functions. Conditions for a
minimum in a most simple way are expressed by means of upper
exhausters and coexhausters, while conditions for a maximum are
described in terms of lower exhausters and coexhausters. Thus the
problem of converting (i.e., obtaining an upper exhauster or
coexhauster when the lower one is given and vice verse) arises.

Numerical examples are provided throughout the paper. Obtained
results can be interesting for researchers in different areas such
as nonsmooth analysis, game theory, mathematical programming and
others.
\end{abstract}

\subsection*{Keywords} Nonsmooth analysis; nondifferentiable optimization;
minmax problems; exhausters; coexhausters.

\section{Introduction}\label{sec1}

To study extremal properties of a function we usually approximate
it with some function from a certain class. Then we derive
optimality conditions and build optimization algorithms in terms
of these approximations. This is a natural way which is used for
smooth functions and can be applied in a nonsmooth case too. Constructive Nonsmooth Analysis is the area based on the ideas introduced by Demyanov and Rubinov. This is one of the many approaches invented for the study of nonsmooth functions.
For example, we can mention Shor
\cite{Demyanov-shor74}, Clarke
\cite{Demyanov-cla83,Clarke_eng}, 
Mordukhovich
\cite{Mordukhovich88,Mordukhovich76,Mordukhovich80,Mordukhovich06,Mordukhovich18},
Michel--Penot \cite{M-P} subdifferentials, approximate and
geometric Ioffe subdifferentials \cite{Io90}, Varga containers
\cite{Warga}. Detailed analysis and more information on these
approaches is given in \cite{BZ}.

In 1980-th Demyanov Rubinov and Polyakova introduced notion of
quasidifferentials
\cite{Demyanov-Rubinov-Polyakova_quasidif79,Demyanov-Rubinov_quasidif80,Demyanov-Vasiliev85}.
These are pairs of convex compact sets which are used to represent
the directional derivative and the approximation of a function at
the point in the form of sum of maximum and minimum of linear
functions. The emergence of quasidifferentials laid the foundation
for constructive nonsmooth analysis. A calculus of
quasidifferentials was developed. The formulas of the calculus
enabled researchers to obtain these pairs for a wide class of
nonsmooth functions. Optimality conditions as well as the
procedures of finding directions of a steepest descent and ascent
when these conditions are not satisfied were formulated in terms
of quasidifferentials \cite{Demyanov-Polyakova80}. Thereafter the
theory of quasidifferentials progressed rapidly due to many
significant studies in the area
\cite{Gorokhovik84,Liqun_Qi90,Pallaschke_Urbanski94,Dempe-Pallaschke97,Stavroulakis09,
Basaeva_Kusraev_Kutateladze16,Dolgopolik_2020,Sukhorukova2017}.

Subsequently exhausters notion appeared as an attempt to expand
the class of studied functions. It was introduced by Demyanov
\cite{Demyanov_optimization_99,Demyanov-dem00a} and is based on
the ideas of Pshenichny \cite{Demyanov-psch80} and Rubinov
\cite{Demyanov-dr82,Demyanov-Rubinov2001}. Lower and upper
exhausters are families of convex compact sets which are used to
describe the directional derivative and the approximation of a
function at the point in the form of maxmin and minmax of linear
functions. Since calculus and optimality conditions have been
derived in terms of exhausters too, this concept retains the
constructiveness of quasidifferentials. At the same time the class
of exhausterable functions is wider than the class of
quasidifferentiable ones. Any quasidifferentiable functions has
exhausters but the opposite is not true.

It turned out that conditions for the minimum most organically are
expressed via an upper exhauster while conditions for a maximum
are described via a lower one. Therefore an upper exhauster is
called proper for the minimization problems and adjoint for the
maximization ones while a lower exhauster is called proper for the
maximization and adjoint for the minimization problems. Therefore
when having adjoint exhauster we can either convert it to get a
proper exhauster or work with the adjoint exhauster itself. The
latter requires that optimality conditions to be derived in terms
of adjoint exhausters. These conditions were first stated by Roshchina \cite{Demyanov-dros05,Demyanov-dros06}. Later this problem
was considered in works of Abbasov
\cite{Demyanov-Abbasov_IMMO10,Demyanov-Abbasov_jogo13,Abbasov_jimo15}
where the conditions were obtained in a
geometrically transparent form, which allows one to get directions
of the steepest descent and ascent by means of adjoint exhausters.

The procedure of exhausters converting  was described in
\cite{Demyanov-dem00a}. It is applicable only to two-dimensional
cases. So the problem of describing more general procedure for
exhausters converting is still opened. In the paper we solve this
problem. The works \cite{Gor_2013,Gor_2015} are also dedicated to the conversion of exhausters.

The fact that, in general, exhauster set-valued mapping is not
continuous in the Hausdorff metric leads to the convergence
problems in algorithms which employ these objects. To overcome
this drawback coexhausters notion was introduced
\cite{Demyanov-dem00a}. These are families of convex compact sets
which describe nonhomogeneous approximations of a function at the
point in the form of maxmin and minmax of affine functions. All
the results obtained for exhausters were generalized and described
for coexhausters, but the problem of converting coexhausters is
open too. One can work with continuous coexhauster set-valued
mapping which guarantee stability and convergence of numerical
algorithms.

The paper is organized as follows: in Section 2 we give
definitions of exhausters and coexhausters and discuss existing
procedure for converting these families; in Section 3 we state and
prove sufficient conditions for equivalence of minmax and maxmin,
also we show that these results do not follow from the basic
minimax theorems; in Section 4 we develop new converting procedure
for exhausters and coexhausters; concluding remarks are presented
in Section 5.

\section{Directional derivative. Exhausters and coexhausters. Converting procedure}

In this section we provide main definitions and results which are
employed in the main part of the paper.

\subsection{Directional derivative}

Let a function $f\colon \mathbb{R}^n \rightarrow  \mathbb{R}$ be
given. The function $f$ is called differentiable at a point $x \in
\mathbb{R}^n$ in a direction $\Delta \in \mathbb{R}^{n}$
if and only if there exists the final limit
\begin{equation*}\label{AM_eq13appr}
f^{\prime}(x,\Delta) = \lim_{\alpha\downarrow 0}\frac{f(x+\alpha
\Delta)-f(x)}{\alpha}.
\end{equation*}

The value $f^{\prime}(x,\Delta)$ is called the directional
derivative of the function $f$ at the point $x\in\mathbb{R}^n$ in
the direction $\Delta\in\mathbb{R}^n$. The directional derivative
$f^{\prime}(x,\Delta)$ is positively homogeneous (p.h.) as
function of direction $\Delta\in\mathbb{R}^n$. This function
represents a p.h. approximation of the increment of the function
$f$ in the neighborhood of the point $x$. Necessary conditions for
a minimum and a maximum are stated in terms of directional
derivatives \cite{rock70}.

\begin{thm}\label{Demyanov_Abbasov_T2.2}
Let a function $f\colon \mathbb{R}^{n}\to \mathbb{R}$ be
directionally differentiable at a point $x_{\ast}\in
\mathbb{R}^n$. For the point $x_{\ast}$ to be a minimizer of the
function $f$ on $\mathbb{R}^{n}$ it is necessary that
\begin{equation}\label{Demyanov_Abbasov_eq16appr}
f^{\prime}(x_{\ast},\Delta)\geq 0, \quad \forall \Delta \in
\mathbb{R}^{n}.
\end{equation}
\end{thm}

\begin{thm}\label{Demyanov_Abbasov_T2.3}
Let a function $f\colon \mathbb{R}^{n}\to \mathbb{R}$ be
directionally differentiable at a point $x^{\ast}\in
\mathbb{R}^n$. For the point $x^{\ast}$ to be a maximizer of the
function $f$ on $\mathbb{R}^{n}$ it is necessary that
\begin{equation}\label{Demyanov_Abbasov_eq20appr}
f^{\prime}(x^{\ast},\Delta)\leq 0, \quad \forall \Delta  \in
\mathbb{R}^{n}.
\end{equation}
\end{thm}

These conditions are not constructive since we cannot check the
sign of directional derivative for all directions. In case of
smooth functions directional derivative can be represented in the
form of a scalar product of the gradient and the direction. Using
this representation optimality conditions can be rewritten in
terms of the gradient. However, this approach is not applicable to
the nonsmooth case since a nonsmooth function cannot be
approximated in the neighborhood of the considered point by a
linear function. So we need another representation of a
directional derivative to obtain constructive optimality
conditions in the nonsmooth case.

\subsection{Exhausters}

Let $f\colon\mathbb{R}^n\to\mathbb{R}$ be a directionally
differentiable function and $h(\Delta)=f^{\prime}(x,\Delta)$ be
the derivative of the function $f$ at a point $x$ in a direction
$\Delta$. Fix $x\in \mathbb{R}^n$. M.Castellani (see \cite{cas98})
proved that, if $h$ is Lipschitz, then there exist families of
convex and compact sets $E^*$ and $E_*$ in the space
$\mathbb{R}^n$ such that $h(\Delta)$ can be written in the form
\begin{equation}\label{eqs3a1}
h(\Delta)=h_1(\Delta)=\min_{C\in E^*} \max_{v\in C}\langle
v,\Delta\rangle, \quad \forall \Delta\in \mathbb{R}^n,
\end{equation} and in the form
\begin{equation}\label{eqs3a2}
h(\Delta)= h_2(\Delta)=\max_{C\in E_*} \min_{w\in C}\langle
w,\Delta\rangle, \quad \forall \Delta\in \mathbb{R}^n.
\end{equation}

The family of sets $E^*$ is called an {\it upper exhauster} of the
function $f$ at the point $x$, while the family  $E_*$ is called a
{\it lower exhauster} of the function $f$ at the point $x$.

For an arbitrary p.h. function $h$ represented in the form
(\ref{eqs3a1}), the family $E^*$ is called an {\it upper
exhauster} of the function $h$. If (\ref{eqs3a2}) holds then the
family $E_*$ is called a {\it lower exhauster} of $h$.

Exhausters were introduced in
\cite{Demyanov_optimization_99,Demyanov-Rubinov2001,Demyanov-dem00a}.
By means of the exhauster representation of the directional
derivative Demyanov reformulated extremal conditions \ref{1904e3a}
and \ref{1904e3b}.

\begin{thm}\label{AM_exhausters_th_unconstrained_min extremum_conditions_upper_ex}
If a function $f(x)$ attains a local minimum at a point $x_\ast$
and an upper exhauster $E^\ast$ of the function $f(x)$ at the
point $x_\ast$ is known, then
$$h(\Delta)=f'(x_\ast,\Delta)=\min_{C\in E^\ast}\max_{v\in C}\langle
v,\Delta\rangle\geq 0 \quad\forall \Delta\in\mathbb{R}^n,$$ what
is equivalent to the condition
\begin{equation}\label{AM_exhausters_th_unconstrained_min_proper_cond}
 0_n\in C \quad \forall \ C\in E^\ast.
\end{equation}
\end{thm}

\begin{thm}\label{AM_exhausters_th_unconstrained_max extremum_conditions}
If a function $f(x)$ attains a local maximum at a point $x^{\ast}$
and a lower exhauster $E_\ast$ of the function $f(x)$ at the point
$x^{\ast}$ is known, then
$$h(\Delta)=f'(x^\ast,\Delta)=\max_{C\in E_\ast}\min_{v\in C}\langle
v,\Delta\rangle\leq 0 \quad\forall g\in\mathbb{R}^n,$$ is
equivalent to the condition
\begin{equation}\label{AM_exhausters_th_unconstrained_max_proper_cond}
 0_n\in C \quad \forall \ C\in E_\ast.
\end{equation}
\end{thm}

Note that the upper and the lower exhausters of the smooth
function $f$ at the point $x$ consist of a singleton which
coincides with the gradient $\nabla f(x)$. Therefore in this case
conditions (\ref{AM_exhausters_th_unconstrained_min_proper_cond})
and (\ref{AM_exhausters_th_unconstrained_max_proper_cond}) means
that the gradient at the considered point equals zero.

Thus, conditions for a minimum are described in terms of upper
exhausters, while conditions for a maximum -- in terms of lower
exhausters. An {\it upper exhauster} is called {\it proper} for
the minimization problem and {\it adjoint} for the maximization
one, while {\it a lower exhauster} is referred to as {\it proper}
for the maximization problem and {\it adjoint } for the
minimization one.

If one studies a minimization problem and only a lower exhauster
is known then it is required to get an upper exhauster. Such an
exhauster can be constructed by means of the converting procedure
which was suggested in \cite{Demyanov-dem00a}.

Let $E\subset 2^{\mathbb{R}^n}$ be a family of convex compacts in
$\mathbb{R}^n$, which is totally bounded, i.e. there exists an
$r<\infty$, such that
$$C\subset B_r(\mathbf{0}) \quad \forall C\in E.$$
For any $\Delta\in\mathbb{R}^n$ such that $\|\Delta\|=1$ build
\begin{equation*}\label{exh_convertors_def}
\widetilde{C}(\Delta)=\operatorname{cl}\operatorname{co}\left\{w(C)\in
C\mid \langle w(C),\Delta\rangle=\min_{w\in C}\langle
w,\Delta\rangle, \ C\in E\right\}.
\end{equation*}
Then the family
$$\widetilde{E}=\left\{C=\widetilde{C}(\Delta) \mid \Delta\in\mathbb{R}^n\right\}$$ is converted family, i.e. if $E$ were a
lower exhauster of the function $h$, then $E^\diamond$ is an upper
exhauster of $h$ and vise versa.

Some problems addressing converting procedure were considered in
\cite{Sang_2017,Demyanov_2011_convert}.

The functions $h_1(x,g)$ and $h_2(x,g)$ in (\ref{eqs3a1}) and
(\ref{eqs3a2}) are discontinuous in the Hausdorff metric as
functions of $x$ and therefore their application in numerical
algorithms brings us to the problems with stability and
convergence. To omit this difficulties we can consider continuous
forms of the representations of the increment. Coexhausters
provide such type of representation. But the loss of the positive
homogeneity of an approximation of the increment is the price we
are paying for this continuity.

\subsection{Coexhausters}

Let a function $f$ be continuous at a point $x\in X$. We say that
at the point $x$ the function $f$ has an upper coexhauster in the
sense if and only if the following expansion holds for any
$\Delta\in\mathbb{R}^n$:
\begin{equation*}\label{1904e3a}
f(x+\Delta)=f(x)+\min_{C\in{\overline{E}(x)}}\max_{[a,v]\in{C}}[a+\langle
v,\Delta\rangle]+o_{x}(\Delta),
\end{equation*}
where $\overline{E}(x)$ is a family of convex compact sets in
${\mathbb{R}}^{n+1}$, and $o_{x}(\Delta)$ satisfies
\begin{equation}\label{exh_Dini_0}\lim_{\alpha\downarrow 0}\frac{o_{x}(\alpha\Delta)}{\alpha}=0, \quad \forall\Delta\in\mathbb{R}^n.\end{equation}

The set $\overline{E}(x)$ is called an upper coexhauster of $f$ at
the point $x$.

We say that at the point $x$ the function $f$ has a lower
coexhauster if and only if the following expansion holds for any
$\Delta\in\mathbb{R}^n$:
\begin{equation*}\label{1904e3b}f(x+\Delta)=f(x)+\max_{C\in{\underline{E}(x)}}\min_{[b,w]\in{C}}[b+\langle w,\Delta\rangle]+o_{x}(\Delta),
\end{equation*}
where $\underline{E}(x)$ is a family of convex compact sets in
${\mathbb{R}}^{n+1}$, and $o_{x}(\Delta)$ satisfies
(\ref{exh_Dini_0}).

The set $\underline{E}(x)$ is called a lower coexhauster of the
function $f$ at the point $x$.

It is obvious that the following equalities holds
$$\min_{C\in{\overline{E}(x)}}\max_{[a,v]\in{C}}a=\max_{C\in{\underline{E}(x)}}\min_{[b,w]\in{C}}b=0,$$
for an upper and a lower coexhauster at any $x$. Therefore we can
deal with the approximations of $f$ itself, i.e.
\begin{equation*}\label{coexh_eqs3a1}
h(\Delta)=h_3(\Delta)=\min_{C\in{\overline{E}}}\max_{[a,v]\in{C}}[a+\langle
v,\Delta\rangle], \quad \forall \Delta\in \mathbb{R}^n,
\end{equation*} and
\begin{equation*}\label{coexh_eqs3a2}
h(\Delta)=
h_4(\Delta)=\max_{C\in{\underline{E}}}\min_{[b,w]\in{C}}[b+\langle
w,\Delta\rangle], \quad \forall \Delta\in \mathbb{R}^n,
\end{equation*}
Hereinafter we consider the families $\overline{E}$ and
$\underline{E}$ as an upper and lower coexhausters of the function
$h$ correspondingly.

The notion of coexhauster was introduced in
\cite{Demyanov_optimization_99,Demyanov-dem00a}, where optimality
conditions were stated in terms of these families.

\begin{thm}\label{T1b}
Let a function $f$ have an upper coexhauster at a point
$x_{\ast}$. For the point ${x}_{\ast}$ to be a minimizer of the
function $f$ on ${\mathbb{R}}^n$, it is necessary that
\begin{equation}\label{1904e3c}
0_{n+1}\in C\quad \forall \ C\in E^{\ast}(x_{\ast}),
\end{equation}
where
\begin{equation}\label{1904e3d}
E^{\ast}(x_{\ast})=\left\{C\in \overline{E}(x_{\ast}) \mid
\max_{[a,v]\in C}a=0\right\}.
\end{equation}
\end{thm}

\begin{thm}\label{T2b}
Let a function $f$ have a lower coexhauster at a point $x^{\ast}$.
For the point $x^{\ast}$ to be a local or global maximizer of the
function $f$ on ${\mathbb{R}}^n$, it is necessary that
\begin{equation}\label{1904e5}
0_{n+1}\in C\quad \forall C\in E_*(x^{\ast}),
\end{equation}
where
\begin{equation}\label{1904e5a}
E_{*}(x^{\ast})=\left\{C\in \underline{E}(x^{\ast})\mid
\min_{[b,w]\in C}b=0\right\}.
\end{equation}
\end{thm}

An upper coexhauster $\overline{E}(x)$ is called a proper one for
the minimization problem (and adjoint for the maximization
problem) while a lower exhauster $\underline{E}(x)$ is called a
proper one for the maximization problem (and adjoint for the
minimization problem). Again we face the problem of obtaining
proper family if adjoint one is given.

If a family $\underline{E}\subset 2^{R^{n+1}}$ of convex compact
sets is a totally bounded lower coexhauster of $h$, then the
family
\begin{equation*}\label{coexh_convertors_E_up_Gamma}
\widetilde{\underline{E}}=\left\{C=C^\diamond(g)\mid
g\in\mathbb{R}^{n+1},\ g_1\geq 0\right\},
\end{equation*}
where
\begin{equation*}\label{coexh_convertors_def_up_C}
C^\diamond(g)=\operatorname{cl}\operatorname{co}\left\{w(C)\in
C\mid \langle w(C),g\rangle=\min_{[b,w]\in C}\langle [b,w],g
\rangle, \ C\in \underline{E}\right\}.
\end{equation*}
is an upper coexhauster of the function $h$.

If a family $\overline{E}\subset 2^{R^n}$ of convex compact sets
is a totally bounded upper coexhauster of $h$, then the family
\begin{equation*}\label{coexh_convertors_E_low_Gamma}
\widetilde{\overline{E}}=\left\{C=C_\diamond(g)\mid
g\in\mathbb{R}^{n+1},\ g_1\geq 0\right\},
\end{equation*}
where
\begin{equation*}\label{coexh_convertors_def_low_C}
C_\diamond(g)=\operatorname{cl}\operatorname{co}\left\{v(C)\in
C\mid \langle v(C),g\rangle=\max_{[a,v]\in C}\langle
[a,v],g\rangle, \ C\in \overline{E}\right\}.
\end{equation*}
is a lower coexhauster of $h$.

For a wide class of functions which have exhausters and
coexhausters these families consists of finite number of convex
polytopes. In what follows we will consider only this case.

The described converting method implies only graphical use.
Construction of converted families is made by means of visual
geometric illustrations. However, this is possible only for
low-dimensional problems. Therefore, the need for more general
conversion procedures arises.

\section{Minimax theorems}

First we state and prove the following result.

\begin{thm}\label{Abbasov_th1} Let $D=\{d_{ij}\}_{k\times p}$ be a
matrix in $\mathbb{R}^{k\times p}$ and there exists
$\overline{j}\in {1,\dots,p}$ such that for all $i\in {1,\dots,k}$
we have
\begin{equation}\label{Abbasov_th1_cond}
d_{i\overline{j}}=\max_{j\in {1,\dots,p}}d_{ij}. \end{equation}
Then the following equation
\begin{equation}\label{Abbasov_th1_main_statemant}\min_{i\in 1,\dots,k}\max_{j\in 1,\dots,p} d_{ij}=\max_{j\in 1,\dots,p}\min_{i\in 1,\dots,k} d_{ij}.\end{equation}
holds.
\end{thm}

\begin{proof}
Denote by $\overline{i}$ the index on which the minimum
$\displaystyle\min_{i\in 1,\dots,k}d_{i\overline{j}}$ is attained,
i.e. $\displaystyle\min_{i\in
1,\dots,k}d_{i\overline{j}}=d_{\overline{i}\overline{j}}$. Then we
have \begin{equation}\label{Abbasov_th1_proof_eq0}\min_{i\in
1,\dots,k}\max_{j\in 1,\dots,p}
d_{ij}=d_{\overline{i}\overline{j}}.\end{equation} Considering
condition (\ref{Abbasov_th1_cond}) we get the chain of
inequalities
$$\min_{i\in 1,\dots,k}d_{ij}\leq d_{\overline{i}j}\leq \max_{j\in 1,\dots,p} d_{\overline{i}j}=d_{\overline{i}\overline{j}}$$
which are true for any $j\in 1,\dots,p$ and therefore we conclude
that
\begin{equation}\label{Abbasov_th1_proof_eq1}
\max_{j\in 1,\dots,p}\min_{i\in 1,\dots,k} d_{ij}\leq
d_{\overline{i}\overline{j}}.
\end{equation}

Since
$\displaystyle\min_{i=1,\dots,k}d_{i\overline{j}}=d_{\overline{i}\overline{j}}$
we obtain
\begin{equation}\label{Abbasov_th1_proof_eq2}
\max_{j\in 1,\dots,p}\min_{i\in 1,\dots,k} d_{ij}\geq
d_{\overline{i}\overline{j}}.
\end{equation}
From (\ref{Abbasov_th1_proof_eq1}) and
(\ref{Abbasov_th1_proof_eq2}) we get that
$$\max_{j\in 1,\dots,p}\min_{i\in 1,\dots,k} d_{ij}=
d_{\overline{i}\overline{j}},$$ whence recalling
(\ref{Abbasov_th1_proof_eq0}) we obtain
(\ref{Abbasov_th1_main_statemant}).
\end{proof}

Similarly can be stated and proved the following theorem.

\begin{thm}\label{Abbasov_th2} Let $D=\{d_{ij}\}_{k\times p}$ be a
matrix in $\mathbb{R}^{k\times p}$ and there exists
$\overline{j}\in {1,\dots,p}$ such that for all $i\in {1,\dots,k}$
we have
\begin{equation}\label{Abbasov_th2_cond}
d_{i\overline{j}}=\min_{j\in {1,\dots,p}}d_{ij}. \end{equation}
Then the following equation
\begin{equation}\label{Abbasov_th2_main_statemant}\max_{i\in 1,\dots,k}\min_{j\in 1,\dots,p} d_{ij}=\min_{j\in 1,\dots,p}\max_{i\in 1,\dots,k} d_{ij}.\end{equation}
holds.
\end{thm}

These two theorems give us sufficient conditions for the equality
of $\min\max$ and $\max\min$ representations. 

Now let us proceed to the more general result which can be considered as a discrete analog of minimax theorem (see \cite{Fan_1953}). 

\begin{thm}\label{FAN} Let $D=\{d_{ij}\}_{k\times p}$ be a
matrix in $\mathbb{R}^{k\times p}$ and there exists
$\widehat{i}\in {1,\dots,k}$ and $\widehat{j}\in {1,\dots,p}$ such
that for all $i\in {1,\dots,k}$ and $j\in {1,\dots,p}$ we have
\begin{equation}\label{FAN_cond}
d_{\widehat{i}j}\leq d_{i\widehat{j}}. \end{equation} Then the
following equation
\begin{equation}\label{main_statemant}\min_{i\in 1,\dots,k}\max_{j\in 1,\dots,p} d_{ij}=\max_{j\in 1,\dots,p}\min_{i\in 1,\dots,k} d_{ij}.\end{equation}
holds.
\end{thm}

\begin{proof}
Condition (\ref{FAN_cond}) implies that
\begin{equation}\label{1}
d_{\widehat{i}j}\leq d_{\widehat{i}\widehat{j}},\quad \forall
j=1,\dots,p \Longleftrightarrow
\max_{j=1,\dots,p}d_{\widehat{i}j}=d_{\widehat{i}\widehat{j}},
\end{equation}
\begin{equation}\label{2}
d_{\widehat{i}\widehat{j}}\leq d_{i\widehat{j}},\quad \forall
i=1,\dots,k \Longleftrightarrow
\min_{i=1,\dots,k}d_{i\widehat{j}}=d_{\widehat{i}\widehat{j}}.
\end{equation}

Therefore (see (\ref{2})) we have $$\max_{j=1,\dots,p}d_{ij}\geq
d_{i\widehat{j}}\geq d_{\widehat{i}\widehat{j}}\quad\forall
i=1,\dots,k,$$ whence
\begin{equation}\label{3}\min_{i=1,\dots,k}\max_{j=1,\dots,p}d_{ij}\geq d_{\widehat{i}\widehat{j}}.\end{equation}

But via (\ref{1}) we can obtain the inequality
\begin{equation}\label{4}\min_{i=1,\dots,k}\max_{j=1,\dots,p}d_{ij}\leq \max_{j=1,\dots,p}d_{\widehat{i}j}=d_{\widehat{i}\widehat{j}}.\end{equation}

From (\ref{3}) and (\ref{4}) we conclude
$$\min_{i=1,\dots,k}\max_{j=1,\dots,p}d_{ij}=d_{\widehat{i}\widehat{j}}.$$

Similarly we can derive the equation

$$\max_{j=1,\dots,p}\min_{i=1,\dots,k}d_{ij}=d_{\widehat{i}\widehat{j}}.$$
\end{proof}

Let us consider two illustrative examples.

\begin{exm}
For the matrix $$D=\begin{pmatrix}
1 & 1 & 3\\
78 & 81 & 100\\
0 & 2 & 7
\end{pmatrix}$$
condition (\ref{Abbasov_th1_cond}) is satisfied
since $\overline{j}=3$ as well as (\ref{FAN_cond}) (we can chose $\widehat{i}=1$,
$\widehat{j}=3$). 
\begin{equation*}\label{Abbasov_th1_main_statemant}\min_{i\in 1,2,3}\max_{j\in 1,2,3} d_{ij}=\max_{j\in 1,2,3}\min_{i\in 1,2,3} d_{ij}=3.\end{equation*}

\end{exm}

\begin{exm}
For the matrix $$D=\begin{pmatrix}
1 & 35 & 15\\
6 & 7 & 10\\
40 & 7 & 20
\end{pmatrix}$$
condition (\ref{FAN_cond}) is satisfied (\,$\widehat{i}=2$,
$\widehat{j}=3$) but (\ref{Abbasov_th1_cond}) is not.
\begin{equation*}\label{Abbasov_th1_main_statemant}\min_{i\in 1,2,3}\max_{j\in 1,2,3} d_{ij}=\max_{j\in 1,2,3}\min_{i\in 1,2,3} d_{ij}=10.\end{equation*}

\end{exm}

\begin{rem}
It is obvious that if condition (\ref{Abbasov_th1_cond}) is satisfied then  (\ref{FAN_cond}) is also true, because we can choose $\widehat{j}=\bar{j}$ and $\widehat{i}$ such that $d_{\widehat{i}\bar{j}}\leq d_{i\bar{j}}$ for all $i$.

However, neither condition (\ref{FAN_cond}) nor the general fact of existence of
the saddle point give us the possibility to construct a conversion
procedure for exhausters and coexhausters. Only conditions of the
form \ref{Abbasov_th1_cond} and \ref{Abbasov_th2_cond} allow us to
do that.
\end{rem}

\section{Converting exhausters and coexhausters}

Theorems \ref{Abbasov_th1} and \ref{Abbasov_th2} can be used in
procedure of converting exhausters and coexhausters. Consider the
problem of obtaining a lower exhauster from an upper one.

\begin{thm}\label{Abbasov_main_th_up_exh}
Let $h\colon \mathbb{R}^n \to \mathbb{R}$ be a function
 such that $\displaystyle h(\Delta)=\min_{C\in E^\ast}\max_{v\in C}\langle
v,\Delta\rangle$ and $E^\ast$ is a finite family of convex compact
sets from $\mathbb{R}^n$, where $E^\ast=\{C_i\mid i=1,\dots,k\}$,
$C_i=co\{v_{ij}\mid j=1,\dots,m_i\}.$ Then the family
$\widetilde{E}^\ast$ which contains $p=m_1m_2\dots m_k$ sets of
the form
$$C=co\{v_{i j_i}\mid i\in 1,\dots,k,\ j_i\in
1,\dots,m_i\}$$ is a lower exhauster of the function $h$.
\end{thm}

\begin{proof} It is obvious that $\widetilde{E}^\ast=\{\widetilde{C}_j\mid
j=1,\dots,p\}$, where $\widetilde{C}_j=co\{\widetilde{v}_{ij}\mid
i\in 1,\dots,k\}$ and $\widetilde{v}_{ij}$ is some vertex of the
set $C_i$.

Choose an arbitrary $\Delta\in\mathbb{R}^n$. Consider the matrix
$D(\Delta)=\{d_{ij}\}_{k\times p}$, which is composed of columns
consisting of inner products of vertices of all of the sets of the
family $\widetilde{E}^\ast$ and $\Delta$, i.e. $d_{ij}=\langle
\widetilde{v}_{ij},\Delta\rangle$. From the way we constructed the
family $\widetilde{E}^\ast$ it is obvious that there exists
$\overline{j}\in {1,\dots,p}$ such that the following condition
holds
$$d_{i\overline{j}}=\max_{j\in {1,\dots,p}}d_{ij}\quad \forall i\in
{1,\dots,k}.$$ Hence due to Theorem \ref{Abbasov_th1} we have
\begin{equation}\label{Abbasov_th3_proof_eq0}\min_{i\in 1,\dots,k}\max_{j\in 1,\dots,p} d_{ij}=\max_{j\in
1,\dots,p}\min_{i\in 1,\dots,k} d_{ij}.\end{equation}

Since
\begin{equation*}\label{Abbasov_th3_proof_eq1}\min_{i\in
1,\dots,k}\max_{j\in 1,\dots,p} d_{ij}=\min_{C\in
E^\ast}\max_{v\in C}\langle v,\Delta\rangle\end{equation*} and
\begin{equation*}\label{Abbasov_th3_proof_eq2}\max_{j\in
1,\dots,k}\min_{i\in 1,\dots,k} d_{ij}=\max_{C\in
\widetilde{E}_\ast}\min_{v\in C}\langle
v,\Delta\rangle\end{equation*} equality
(\ref{Abbasov_th3_proof_eq0}) implies that the family
$\widetilde{E}^\ast$ is a lower exhauster of the function $h$.
\end{proof}

Similarly we can state and prove the following results for
converting a lower exhauster and upper and lower coexhausters.

\begin{thm}\label{Abbasov_main_th_low_exh}
Let $h\colon \mathbb{R}^n \to \mathbb{R}$ be a function
 such that $\displaystyle h(\Delta)=\max_{C\in E_\ast}\min_{w\in C}\langle
w,\Delta\rangle$ and $E_\ast$ is a finite family of convex compact
sets from $\mathbb{R}^n$, where $E_\ast=\{C_i\mid i=1,\dots,k\}$,
$C_i=co\{w_{ij}\mid j=1,\dots,m_i\}.$ Then the family
$\widetilde{E}_\ast$ which contains $p=m_1m_2\dots m_k$ sets of
the form
$$C=co\{w_{i j_i}\mid i\in 1,\dots,k,\ j_i\in
1,\dots,m_i\}$$ is an upper exhauster of the function $h$.
\end{thm}

\begin{thm}\label{Abbasov_main_th_up_coexh}
Let $h\colon \mathbb{R}^n \to \mathbb{R}$ be a function
 such that $\displaystyle h(\Delta)=\min_{C\in \overline{E}}\max_{[a,v]\in C}[a+\langle
v,\Delta\rangle]$ and $\overline{E}$ is a finite family of convex
compact sets from $\mathbb{R}^{n+1}$, where
$\overline{E}=\{C_i\mid i=1,\dots,k\}$,
$C_i=co\{[a_{ij},v_{ij}]\mid j=1,\dots,m_i\}.$ Then the family
$\widetilde{\overline{E}}$ which contains $p=m_1m_2\dots m_k$ sets
of the form
$$C=co\{[a_{i j_i}, v_{i j_i}]\mid i\in 1,\dots,k,\ j_i\in
1,\dots,m_i\}$$ is a lower coexhauster of the function $h$.
\end{thm}

\begin{thm}\label{Abbasov_main_th_low_coexh}
Let $h\colon \mathbb{R}^n \to \mathbb{R}$ be a function such that
$\displaystyle h(\Delta)=\max_{C\in \underline{E}}\min_{[b,w]\in
C}[b+\langle w,\Delta\rangle]$ and $\underline{E}$ is a finite
family of convex compact sets from $\mathbb{R}^{n+1}$, where
$\underline{E}=\{C_i\mid i=1,\dots,k\}$,
$C_i=co\{[b_{ij},w_{ij}]\mid j=1,\dots,m_i\}.$ Then the family
$\widetilde{\underline{E}}$ which contains $p=m_1m_2\dots m_k$
sets of the form
$$C=co\{[b_{i j_i}, w_{i j_i}]\mid i\in 1,\dots,k,\ j_i\in
1,\dots,m_i\}$$ is an upper coexhauster of the function $h$.
\end{thm}

Let us consider some illustrative examples.

\begin{exm} Let the function $h\colon\mathbb{R}^4\to\mathbb{R}$, such
that
$$h(\Delta)=\max_{C\in E_\ast}\min_{w\in C}\langle
w,\Delta\rangle$$ be given, where $E=\{C_1,C_2\}$ and
$$C_1=\left\{\begin{pmatrix}-&1\\&1\\&1\\&1\end{pmatrix},\begin{pmatrix}&1\\&1\\&1\\&1\end{pmatrix}\right\},\quad C_2=\left\{\begin{pmatrix}&1\\-&1\\-&1\\-&1\end{pmatrix},\begin{pmatrix}-&1\\-&1\\-&1\\-&1\end{pmatrix}\right\}.$$
Using Theorem \ref{Abbasov_main_th_low_exh} we get an upper
exhauster of the form
$\widetilde{E}_\ast=\{\widetilde{C}_1,\widetilde{C}_2,\widetilde{C}_3,\widetilde{C}_4\}$,
$$\widetilde{C}_1=\left\{\begin{pmatrix}-&1\\&1\\&1\\&1\end{pmatrix},\begin{pmatrix}&1\\-&1\\-&1\\-&1\end{pmatrix}\right\},\quad
\widetilde{C}_2=\left\{\begin{pmatrix}-&1\\&1\\&1\\&1\end{pmatrix},\begin{pmatrix}-&1\\-&1\\-&1\\-&1\end{pmatrix}\right\},$$
$$\widetilde{C}_3=\left\{\begin{pmatrix}&1\\&1\\&1\\&1\end{pmatrix},\begin{pmatrix}&1\\-&1\\-&1\\-&1\end{pmatrix}\right\},\quad
\widetilde{C}_4=\left\{\begin{pmatrix}&1\\&1\\&1\\&1\end{pmatrix},\begin{pmatrix}-&1\\-&1\\-&1\\-&1\end{pmatrix}\right\}.$$
Thus $h$ can be represented as
$$h(\Delta)=\min_{C\in \widetilde{E}_\ast}\max_{v\in C}\langle
v,\Delta\rangle.$$ \end{exm}

\begin{exm}
Consider the function $h\colon\mathbb{R}^4\to\mathbb{R}$,
$$h(\Delta)=\min_{C\in \overline{E}}\max_{[a,v]\in C}[a+\langle
v,\Delta\rangle],$$ where $\overline{E}=\{C_1,C_2\}$
$$C_1=\operatorname{co}\{[a_i,v_i]\mid i=1,2,3\},\quad C_2=\{\mathbf{0}_4\},$$
$a_i=1$ for all $i=1,2,3,$ and $v_i$ is $i$-th standard basis
vector, i.e. $v_i=e_i$ for all $i=1,2,3$.

Via Theorem \ref{Abbasov_main_th_up_coexh} we obtain a lower
coexhauster of the form
$\widetilde{\overline{E}}=\{\widetilde{C}_1,\widetilde{C}_2,\widetilde{C}_3\}$,
where
$$\widetilde{C}_i=\operatorname{co}\{[a_i,e_i],\mathbf{0}_4\}\quad\forall i=1,2,3,$$
and therefore can represent the function $h$ as
$$h(\Delta)=\max_{C\in \widetilde{\overline{E}}}\min_{[b,w]\in C}[b+\langle
v,\Delta\rangle].$$

\end{exm}

\section{Conclusion} Obtained results give a solution of the
problem of converting exhausters and coexhausters in cases when
these families consist of finite number of convex polytopes. 
After applying the developed procedure we usually get a family
with many redundant sets. These sets can be discarded via various
reduction techniques and methods presented in
\cite{Grzybowski_Pallaschke_2010,Roshchina_2008_JOTA_reduction,Roshchina_2008_JCA_reduction,Abbasov_JOGO_2019}.

\end{document}